\documentclass{amsart}

\usepackage[utf8]{inputenc}
\usepackage[english]{babel}
\usepackage[T1]{fontenc}

\usepackage{hyperref}

\usepackage{amsmath}
\usepackage{amsfonts}
\usepackage{amsthm}

\theoremstyle{plain}

\newtheorem{theorem}{Theorem}
\newtheorem{lemma}[theorem]{Lemma}

\theoremstyle{remark}

\newtheorem*{remark}{Remark}

\usepackage{theoremref}
\usepackage{colonequals}
\usepackage{cancel}

\begin{document}

\title[Hessian determinant of quadratic hypersurface]{A concise formula for the Hessian determinant of a~function parameterising a~quadratic hypersurface}

\author[B. Zawalski]{Bartłomiej Zawalski}
\address{Polish Academy of Sciences, Institute of Mathematics, ul. Jana i Jędrzeja Śniadeckich 8, 00-656 Warszawa, Poland}
\email{b.zawalski@impan.pl}
\thanks{This work was supported by the Polish National Science Centre grants no. 2015/18/A/ST1/00553 and 2020/02/Y/ST1/00072. Numerical computations were performed within the Interdisciplinary Centre for Mathematical and Computational Modelling UW grant number G76-28.}

\begin{abstract}
We will give a~concise formula for the Hessian determinant of a~smooth function $y:\mathbb R^n\supseteq\Omega\to\mathbb R$ such that its graph is contained in a~quadratic hypersurface. The proof will make heavy use of matrix algebra.
\end{abstract}

\maketitle

The aim of this work is to prove the following theorem:

\begin{theorem}\thlabel{thm:01}
Let $y:\mathbb R^n\supseteq\Omega\to\mathbb R$ be a~function of class $C^2$ defined on an open subset of $\mathbb R^n$ and satisfying a~quadratic equation
\begin{equation}\label{eq:01}\mathbf v(\mathbf x)^\top\mathbf Q\mathbf v(\mathbf x)=0,\end{equation}
where
$$\mathbf v(\mathbf x)\colonequals\begin{bmatrix}\mathbf x \\ y(\mathbf x) \\ 1\end{bmatrix},\quad \mathbf x\in\Omega$$
is an augmented $(n+2)$-dimensional column vector and $\mathbf Q$ is an arbitrary $\left(n+2\right)\times\left(n+2\right)$ matrix. Then the following formula holds:
\begin{equation}\label{eq:02}|\mp\mathbf H_y(\mathbf x)|\cdot\Delta_y(\mathbf x)^{n/2+1}=-\left|\mathbf Q+\mathbf Q^\top\right|,\end{equation}
where $\mathbf H_y$ is the Hessian matrix of $y$, $\Delta_y$ is the discriminant of the left hand side of \eqref{eq:01} with respect to the variable $y$ and the sign $\mp$ depends on the selected branch of the square root function.
\end{theorem}

Here by \emph{Hessian matrix} of $y$ we mean a~square $n\times n$ matrix defined as follows:
\begin{equation}\label{eq:13}\mathbf H_y\colonequals\begin{bmatrix}
  \dfrac{\partial^2 y}{\partial x_1^2} & \dfrac{\partial^2 y}{\partial x_1\,\partial x_2} & \cdots & \dfrac{\partial^2 y}{\partial x_1\,\partial x_n} \\[2.2ex]
  \dfrac{\partial^2 y}{\partial x_2\,\partial x_1} & \dfrac{\partial^2 y}{\partial x_2^2} & \cdots & \dfrac{\partial^2 y}{\partial x_2\,\partial x_n} \\[2.2ex]
  \vdots & \vdots & \ddots & \vdots \\[2.2ex]
  \dfrac{\partial^2 y}{\partial x_n\,\partial x_1} & \dfrac{\partial^2 y}{\partial x_n\,\partial x_2} & \cdots & \dfrac{\partial^2 y}{\partial x_n^2}
\end{bmatrix},\end{equation}
and by \emph{discriminant} of the quadratic polynomial $ay^2+by+c$ with $a\neq 0$ with respect to the variable $y$ we mean
$$\Delta_y\colonequals b^2-4ac.$$

In the course of the proof we will also use a~vector differential operator, usually represented by the nabla symbol $\nabla$. It is defined in terms of partial derivative operators as
$$\nabla=\begin{bmatrix}\frac{\partial}{\partial x_1}\\\vdots\\\frac{\partial}{\partial x_n}\end{bmatrix}.$$
In a~convenient mathematical notation, we can consider e.g. a~formal product of nabla with a~scalar and a~formal tensor product of nabla with a~vector field. Namely,
$$\nabla y\colonequals\begin{bmatrix}\frac{\partial y}{\partial x_1}\\\vdots\\\frac{\partial y}{\partial x_n}\end{bmatrix},\quad\nabla\otimes\vec v\colonequals\begin{bmatrix}\frac{\partial\vec v_1}{\partial x_1}&\cdots&\frac{\partial\vec v_n}{\partial x_1}\\\vdots&\ddots&\vdots\\\frac{\partial\vec v_1}{\partial x_n}&\cdots&\frac{\partial\vec v_n}{\partial x_n}\end{bmatrix}.$$
In particular,
$$\mathbf H_y=\nabla\otimes\nabla y.$$

\begin{proof}
To simplify the notation, we will henceforth skip the dependence on $\mathbf x$ whenever it is clear from the context.\\

Since the formula \eqref{eq:02} depends on $\mathbf Q$ only through $\mathbf Q+\mathbf Q^\top$, without loss of generality we may assume that $\mathbf Q$ is symmetric by taking its symmetric part
$$\frac{1}{2}\left(\mathbf Q+\mathbf Q^\top\right).$$
Thus we can write
$$\mathbf Q\equalscolon\begin{bmatrix} \mathbf A=\mathbf A^\top & \mathbf b & \mathbf c \\ \mathbf b^\top & d & e \\ \mathbf c^\top & e & f \end{bmatrix}$$
as a~symmetric block matrix, where $\mathbf A$ is a~symmetric $n\times n$ matrix, $\mathbf b,\mathbf c$ are $n$-dimensional column vectors and $d,e,f$ are scalars.\\

In terms of these new variables, equation \eqref{eq:01} reads
\begin{equation}\label{eq:03}\begin{bmatrix} \mathbf x \\ y \\ 1 \end{bmatrix}^\top
\begin{bmatrix} \mathbf A & \mathbf b & \mathbf c \\ \mathbf b^\top & d & e \\ \mathbf c^\top & e & f \end{bmatrix}
\begin{bmatrix} \mathbf x \\ y \\ 1 \end{bmatrix}
=dy^2+2\left(\mathbf b^\top\mathbf x+e\right)y+\left(\mathbf x^\top\mathbf A\mathbf x+2\mathbf c^\top\mathbf x+f\right)\end{equation}
and hence its discriminant with respect to the variable $y$ is itself a~quadratic polynomial in $\mathbf x$ given by
\begin{align}
\label{eq:14}\Delta_y&=4\left(\mathbf b^\top\mathbf x+e\right)^2-4d\left(\mathbf x^\top\mathbf A\mathbf x+2\mathbf c^\top\mathbf x+f\right)\\
\nonumber&=4\mathbf x^\top\left(\mathbf b\mathbf b^\top-d\mathbf A\right)\mathbf x+8\left(e\mathbf b^\top-d\mathbf c^\top\right)\mathbf x+4\left(e^2-df\right).
\end{align}
Denote its coefficients by
\begin{equation}\label{eq:11}\boldsymbol\Lambda\colonequals d\mathbf A-\mathbf b\mathbf b^\top,\quad\boldsymbol\mu\colonequals e\mathbf b-d\mathbf c,\quad\nu\colonequals df-e^2\end{equation}
respectively, so that the equality
\begin{equation}\label{eq:04}\Delta_y\equalscolon-4\mathbf x^\top\boldsymbol\Lambda\mathbf x+8\boldsymbol\mu^\top\mathbf x-4\nu \end{equation}
holds.\\

Since the formula \eqref{eq:02} is continuous in both $\mathbf x$ and $\mathbf Q$, it is enough to prove it for a~dense subset of pairs $(\mathbf x,\mathbf Q)$. Hence without loss of generality we may assume that $d\neq 0$ and $\Delta_y>0$. Solving \eqref{eq:03} for $y$ yields
$$y=\frac{-2\left(\mathbf b^\top\mathbf x+e\right)\pm\sqrt{\Delta_y}}{2d},$$
where the sign $\pm$ depends on the selected branch of the square root function and is opposite to that in \eqref{eq:02}. Since the Hessian matrix does not depend on the linear part, we have
$$\pm\mathbf H_y=\frac{\nabla\otimes\nabla\sqrt{\Delta_y}}{2d}$$
with
\begin{equation}\label{eq:08}|\pm\mathbf H_y|=2^{-n}d^{-n}\left|\nabla\otimes\nabla\sqrt{\Delta_y}\right|.\end{equation}

Now, for an arbitrary function $\Delta_y$, the Hessian matrix of $\sqrt{\Delta_y}$ is given by
$$\nabla\otimes\nabla\sqrt{\Delta_y}=\nabla\otimes\frac{\nabla\Delta_y}{2{\Delta_y}^{1/2}}=\frac{2\Delta_y\nabla\otimes\nabla\Delta_y-\nabla\Delta_y\otimes\nabla\Delta_y}{4{\Delta_y}^{3/2}}$$
and thus is a~linear combination of a~square matrix $\nabla\otimes\nabla\Delta_y$ and a~rank one matrix $\nabla\Delta_y\otimes\nabla\Delta_y$. Its determinant can be computed using the following simple fact from linear algebra:

\begin{lemma}[{\cite[Theorem~18.1.1]{harville2011matrix}}]\thlabel{lem:01}
Let $\mathbf R$ represent an $n\times n$ matrix, $\mathbf S$ an $n\times m$ matrix, $\mathbf T$ an $m\times m$ matrix, and $\mathbf U$ an $m\times n$ matrix. If $\mathbf R$ and $\mathbf T$ are non-singular, then
$$|\mathbf R+\mathbf{STU}|=|\mathbf R||\mathbf T||\mathbf T^{-1}+\mathbf U\mathbf R^{-1}\mathbf S|.$$
\end{lemma}

Recall that $\Delta_y$ is a~quadratic polynomial \eqref{eq:04}, in which case
$$\nabla\Delta_y=-8\boldsymbol\Lambda\mathbf x+8\boldsymbol\mu,\quad\nabla\otimes\nabla\Delta_y=-8\boldsymbol\Lambda.$$
Again, without loss of generality we may assume that $\boldsymbol\Lambda$ is non-singular. Applying \thref{lem:01} to
$$\mathbf R\colonequals 2\Delta_y\nabla\otimes\nabla\Delta_y,\quad\mathbf S\colonequals-(\nabla\Delta_y),\quad\mathbf T\colonequals\mathbf I_1,\quad\mathbf U\colonequals(\nabla\Delta_y)^\top$$
yields
\begin{equation}\label{eq:09}
\begin{aligned}
&\left|\nabla\otimes\nabla\sqrt{\Delta_y}\right|\\
&\quad=4^{-n}{\Delta_y}^{-3n/2}|\mathbf R|\left(1-(\nabla\Delta_y)^\top\mathbf R^{-1}(\nabla\Delta_y)\right)\\
&\quad=2^{-n-1}{\Delta_y}^{-n/2-1}|\nabla\otimes\nabla\Delta_y|\left(2\Delta_y-(\nabla\Delta_y)^\top(\nabla\otimes\nabla\Delta_y)^{-1}(\nabla\Delta_y)\right).
\end{aligned}
\end{equation}
Further, without loss of generality we may assume that $\mathbf A$ is non-singular. Applying \thref{lem:01} to
$$\mathbf R\colonequals d\mathbf A,\quad\mathbf S\colonequals-\mathbf b,\quad\mathbf T\colonequals\mathbf I_1,\quad\mathbf U\colonequals\mathbf b^\top$$
yields
\begin{equation}\label{eq:05}
\begin{aligned}
|\nabla\otimes\nabla\Delta_y|&=(-8)^n|\mathbf R|\left(1-\mathbf b^\top\mathbf R^{-1}\mathbf b\right)\\
&=(-8)^nd^{n-1}|\mathbf A|\left(d-\mathbf b^\top\mathbf A^{-1}\mathbf b\right).
\end{aligned}
\end{equation}
Moreover,
$$(\nabla\Delta_y)^\top(\nabla\otimes\nabla\Delta_y)^{-1}(\nabla\Delta_y)=-8\mathbf x^\top\boldsymbol\Lambda\mathbf x+16\boldsymbol\mu^\top\mathbf x-8\boldsymbol\mu^\top\boldsymbol\Lambda^{-1}\boldsymbol\mu$$
and consequently
\begin{equation}\label{eq:06}2\Delta_y-(\nabla\Delta_y)^\top(\nabla\otimes\nabla\Delta_y)^{-1}(\nabla\Delta_y)=-8\nu+8\boldsymbol\mu^\top\boldsymbol\Lambda^{-1}\boldsymbol\mu \end{equation}
is a~constant independent of $\mathbf x$. Since $\boldsymbol\Lambda=d\mathbf A-\mathbf b\mathbf b^\top$ is a~linear combination of a~square matrix $\mathbf A$ and a~rank one matrix $\mathbf b\mathbf b^\top$, its inverse can be computed using another simple fact from linear algebra:

\begin{lemma}[{\cite[Corollary~18.2.10]{harville2011matrix}}]\thlabel{lem:02}
Let $\mathbf R$ represent an $n\times n$ non-singular matrix, and let $\mathbf s$ and $\mathbf u$ represent $n$-dimensional column vectors. Then $\mathbf R+\mathbf s\mathbf u^\top$ is non-singular if and only if $\mathbf u^\top\mathbf R^{-1}\mathbf s\neq-1$, in which case
$$\left(\mathbf R+\mathbf s\mathbf u^\top\right)^{-1}=\mathbf R^{-1}-\left(1+\mathbf u^\top\mathbf R^{-1}\mathbf s\right)^{-1}\mathbf R^{-1}\mathbf s\mathbf u^\top\mathbf R^{-1}.$$
\end{lemma}

Applying \thref{lem:02} to
$$\mathbf R\colonequals d\mathbf A,\quad\mathbf s\colonequals-\mathbf b,\quad\mathbf u\colonequals\mathbf b$$
yields
\begin{equation}\label{eq:07}
\begin{aligned}
\boldsymbol\Lambda^{-1}&=\mathbf R^{-1}+\left(1-\mathbf b^\top\mathbf R^{-1}\mathbf b\right)^{-1}\mathbf R^{-1}\mathbf b\mathbf b^\top\mathbf R^{-1}\\
&=d^{-1}\left(\mathbf A^{-1}+\left(d-\mathbf b^\top\mathbf A^{-1}\mathbf b\right)^{-1}\mathbf A^{-1}\mathbf b\mathbf b^\top\mathbf A^{-1}\right)\\
&=d^{-1}\left(d-\mathbf b^\top\mathbf A^{-1}\mathbf b\right)^{-1}\left(\left(d-\mathbf b^\top\mathbf A^{-1}\mathbf b\right)\mathbf A^{-1}+\mathbf A^{-1}\mathbf b\mathbf b^\top\mathbf A^{-1}\right).
\end{aligned}
\end{equation}
By combining \eqref{eq:05} and \eqref{eq:06} we obtain
\begin{align*}
&|\nabla\otimes\nabla\Delta_y|\left(2\Delta_y-(\nabla\Delta_y)^\top(\nabla\otimes\nabla\Delta_y)^{-1}(\nabla\Delta_y)\right)\\
&\quad=(-8)^nd^{n-1}|\mathbf A|\left(d-\mathbf b^\top\mathbf A^{-1}\mathbf b\right)\left(-8\nu+8\boldsymbol\mu^\top\boldsymbol\Lambda^{-1}\boldsymbol\mu\right)\\
&\quad=(-8)^{n+1}d^{n-2}|\mathbf A|\cdot d\left(d-\mathbf b^\top\mathbf A^{-1}\mathbf b\right)\left(\nu-\boldsymbol\mu^\top\boldsymbol\Lambda^{-1}\boldsymbol\mu\right).
\end{align*}
Further, using \eqref{eq:07} gives us
\begin{align*}
&d\left(d-\mathbf b^\top\mathbf A^{-1}\mathbf b\right)\left(\nu-\boldsymbol\mu^\top\boldsymbol\Lambda^{-1}\boldsymbol\mu\right)\\
&\quad=d\left(d-\mathbf b^\top\mathbf A^{-1}\mathbf b\right)\nu-\boldsymbol\mu^\top\left(d\left(d-\mathbf b^\top\mathbf A^{-1}\mathbf b\right)\boldsymbol\Lambda^{-1}\right)\boldsymbol\mu\\
&\quad=d\left(d-\mathbf b^\top\mathbf A^{-1}\mathbf b\right)\nu-\boldsymbol\mu^\top\left(\left(d-\mathbf b^\top\mathbf A^{-1}\mathbf b\right)\mathbf A^{-1}+\mathbf A^{-1}\mathbf b\mathbf b^\top\mathbf A^{-1}\right)\boldsymbol\mu,
\end{align*}
which after applying the definitions \eqref{eq:11} expands to
\begin{align*}
&d^3f-d^2e^2-d^2f\left(\mathbf b^\top\mathbf A^{-1}\mathbf b\right)+\cancel{de^2\left(\mathbf b^\top\mathbf A^{-1}\mathbf b\right)}-\cancel{de^2\left(\mathbf b^\top\mathbf A^{-1}\mathbf b\right)}\\
&\quad+\cancel{e^2\left(\mathbf b^\top\mathbf A^{-1}\mathbf b\right)\left(\mathbf b^\top\mathbf A^{-1}\mathbf b\right)}-\cancel{e^2\left(\mathbf b^\top\mathbf A^{-1}\mathbf b\right)\left(\mathbf b^\top\mathbf A^{-1}\mathbf b\right)}+d^2e\left(\mathbf c^\top\mathbf A^{-1}\mathbf b\right)\\
&\quad-\cancel{de\left(\mathbf b^\top\mathbf A^{-1}\mathbf b\right)\left(\mathbf c^\top\mathbf A^{-1}\mathbf b\right)}+\cancel{de\left(\mathbf c^\top\mathbf A^{-1}\mathbf b\right)\left(\mathbf b^\top\mathbf A^{-1}\mathbf b\right)}+d^2e\left(\mathbf b^\top\mathbf A^{-1}\mathbf c\right)\\
&\quad-\cancel{de\left(\mathbf b^\top\mathbf A^{-1}\mathbf b\right)\left(\mathbf b^\top\mathbf A^{-1}\mathbf c\right)}+\cancel{de\left(\mathbf b^\top\mathbf A^{-1}\mathbf b\right)\left(\mathbf b^\top\mathbf A^{-1}\mathbf c\right)}-d^3\left(\mathbf c^\top\mathbf A^{-1}\mathbf c\right)\\
&\quad+d^2\left(\mathbf b^\top\mathbf A^{-1}\mathbf b\right)\left(\mathbf c^\top\mathbf A^{-1}\mathbf c\right)-d^2\left(\mathbf c^\top\mathbf A^{-1}\mathbf b\right)\left(\mathbf b^\top\mathbf A^{-1}\mathbf c\right)
\end{align*}
and then after cancellation can be put in the form $d^2\xi$, where
$$\xi\colonequals\begin{vmatrix} d-\mathbf b^\top\mathbf A^{-1}\mathbf b & e-\mathbf b^\top\mathbf A^{-1}\mathbf c \\ e-\mathbf c^\top\mathbf A^{-1}\mathbf b & f-\mathbf c^\top\mathbf A^{-1}\mathbf c \end{vmatrix}.$$
Thus we have eventually arrived at
\begin{equation}\label{eq:12}
\begin{aligned}
&|\nabla\otimes\nabla\Delta_y|\left(2\Delta_y-(\nabla\Delta_y)^\top(\nabla\otimes\nabla\Delta_y)^{-1}(\nabla\Delta_y)\right)\\
&\quad=(-8)^{n+1}d^{n-2}|\mathbf A|\cdot d^2\xi\\
&\quad=(-8)^{n+1}d^n|\mathbf A|\xi.
\end{aligned}
\end{equation}

On the other hand, since we assumed $\mathbf A$ to be non-singular, by means of a~finite sequence of elementary column operations we are able to eliminate $\mathbf b$ and $\mathbf c$ from the last two columns of $\mathbf Q$. Indeed, both
$$\begin{bmatrix} \mathbf A \\ \mathbf b^\top \\ \mathbf c^\top \end{bmatrix}\mathbf A^{-1}\mathbf b=\begin{bmatrix} \mathbf b \\ \mathbf b^\top\mathbf A^{-1}\mathbf b \\ \mathbf c^\top\mathbf A^{-1}\mathbf b \end{bmatrix}$$
and
$$\begin{bmatrix} \mathbf A \\ \mathbf b^\top \\ \mathbf c^\top \end{bmatrix}\mathbf A^{-1}\mathbf c=\begin{bmatrix} \mathbf c \\ \mathbf b^\top\mathbf A^{-1}\mathbf c \\ \mathbf c^\top\mathbf A^{-1}\mathbf c \end{bmatrix}$$
are linear combinations of leading $n$ columns of $\mathbf Q$ with vectors of coefficients $\mathbf A^{-1}\mathbf b$ and $\mathbf A^{-1}\mathbf c$, respectively. Now, since the determinant is invariant under column addition, it follows that
$$|\mathbf Q|=\begin{vmatrix} \mathbf A & \mathbf b & \mathbf c \\ \mathbf b^\top & d & e \\ \mathbf c^\top & e & f \end{vmatrix}=\begin{vmatrix} \mathbf A & \mathbf 0 & \mathbf 0 \\ \mathbf b^\top & d-\mathbf b^\top\mathbf A^{-1}\mathbf b & e-\mathbf b^\top\mathbf A^{-1}\mathbf c \\ \mathbf c^\top & e-\mathbf c^\top\mathbf A^{-1}\mathbf b & f-\mathbf c^\top\mathbf A^{-1}\mathbf c \end{vmatrix}.$$
Observe that the latter matrix is block-lower-triangular, in which case its determinant is simply equal to $|\mathbf A|\xi$. Therefore \eqref{eq:12} reads
\begin{equation}\label{eq:10}|\nabla\otimes\nabla\Delta_y|\left(2\Delta_y-(\nabla\Delta_y)^\top(\nabla\otimes\nabla\Delta_y)^{-1}(\nabla\Delta_y)\right)=(-8)^{n+1}d^n|\mathbf Q|.\end{equation}

Finally, combining \eqref{eq:08}, \eqref{eq:09} and \eqref{eq:10} yields
\begin{align*}
|\pm\mathbf H_y|&=2^{-n}d^{-n}\cdot 2^{-n-1}{\Delta_y}^{-n/2-1}\cdot (-8)^{n+1}d^n|\mathbf Q|\\
&=(-1)^n{\Delta_y}^{-n/2-1}\cdot-2^{n+2}|\mathbf Q|,
\end{align*}
which for symmetric matrix $\mathbf Q$ is equivalent to \eqref{eq:02}. This concludes the proof.
\end{proof}

\begin{remark}
Observe that if $|\mathbf H_y|$ or $\Delta_y$ takes a~zero value anywhere on $\Omega$, then the matrix $\mathbf Q$ is singular and thus $|\mathbf H_y|$ vanishes identically. Indeed, if $\Delta_y$ vanishes identically, then $\mathbf Q$ is a~rank-one matrix, in which case $y$ is affine and thus its Hessian determinant is zero. Othwerwise, as a~non-zero quadratic polynomial, $\Delta_y$ does not vanish on an open dense subset of $\Omega$, which implies that $|\mathbf H_y|$ must be zero there. Hence, by continuity, it vanishes on the whole $\Omega$.
\end{remark}

\begin{remark}
The sign $\mp$ which depends on the selected branch of the square root function has any significance only for odd $n$.
\end{remark}

\begin{remark}
Since the proof of \thref{thm:01} was purely algebraic, the same result holds if we consider e.g. the field $\mathbb C$ of complex numbers instead of the field $\mathbb R$ of real numbers. However, note that both \eqref{eq:01} and \eqref{eq:13} should then be interpreted in the holomorphic sense, i.e. without conjugation.
\end{remark}

\bibliography{references}{}
\bibliographystyle{amsplain}

\end{document}